\documentclass{tran-l}

\def\sub{\substack}

\def\bi{\bibitem}
\def\lb{\langle}
\def\rb{\rangle}

\newcommand{\al}{\alpha}
\def\be{\beta}
\def\sph{\mathbb{S}^{d-1}}
\def \b{\beta}

\def\bl{\bigl}
\def\br{\bigr}

\def\dmin{\displaystyle \min}
\def\Ld{\Lambda}
\def\da{\delta}

\def\df{\displaystyle\frac}

\def\dsum{\displaystyle\sum}
\def\dint{\displaystyle\int}

\def\sub{\substack}
\def\Bl{\Bigl}
\def\Br{\Bigr}
\def\f{\frac}
\def\({\Bigl(}
\def \){ \Bigr)}

\def\ga{\gamma}
\def\Ga{\Gamma}

\def\hb{\hfill$\Box$}

\newcommand{\ld}{\lambda}

\def\sa{\sigma}
\def\dmax{\displaystyle\max}

\def\Og{\Omega}

\def\ta{\theta}

\def\vi{\varphi}

\def\wt{\widetilde}

\newtheorem{thm}{Theorem}[section]

\newtheorem{lem}[thm]{Lemma}

\theoremstyle{definition}
\theoremstyle{remark}
\numberwithin{equation}{section}

\theoremstyle{remark}

\begin{document}

\title{Characterizations of function  spaces  on  the sphere using   frames   }
\author{Feng Dai}
\address{Department of Mathematical and Statistical Sciences, CAB 632,
University of Alberta, Edmonton, Alberta, T6G 2G1,  Canada.}
\email{dfeng@math.ualberta.ca}
\thanks{The  author was supported in part by the NSERC  Canada under grant G121211001.}

\begin{abstract}
In this paper we introduce a polynomial  frame  on the unit sphere
$\sph$ of $\mathbb{R}^d$,  for which every distribution  has a
wavelet-type decomposition.  More importantly, we prove that many
function spaces on the sphere $\sph$, such as $L^p$, $H^p$ and
Besov spaces, can be characterized  in terms of the coefficients
in the wavelet
 decompositions, as in the usual Euclidean case
$\mathbb{R}^d$.  We also study a related nonlinear $m$-term
approximation problem on $\sph$. In particular,  we prove both a
Jackson--type inequality and a   Bernstein--type inequality
associated to  wavelet decompositions, which extend the
corresponding results obtained by R. A. DeVore, B. Jawerth and V.
Popov (``Compression of wavelet decompositions'', {\it  Amer. J.
Math. } {\bf 114} (1992), no. 4, 737--785).
\end{abstract}
 \maketitle

\section{Introduction and summary of main  results}

\subsection{Notations and basic facts.}
We start with some necessary notations.    Given an integer $d\ge
3$, we denote by  $\sph$   the unit sphere of the $d$-dimensional
Euclidean space $\mathbb{R}^d$ and  $d\sa(x)$ the usual Lebesgue
measure on $\sph$ normalized by $\int_{\sph} \,d\sa(x)=1$. For
$0<p\leq\infty$, we let  $L^p\equiv L^p(\sph)$ denote  the usual
Lebesgue space on $\sph$ endowed with the quasi-norm $\|\cdot\|_p$
and $H^p\equiv H^p(\sph)$  the usual Hardy space on $\sph$ endowed
with the quasi-norm $\|\cdot\|_{H^p}$ ( see Section 4 for precise
definition of $H^p(\sph)$). Given a
 measurable subset $E$ of $\sph$,
 we denote by   $|E|$   its  Lebesgue
measure  and  $\chi_E$  its  characteristic function. For $x,
y\in\sph$, we shall use the notation  $d(x,y)$ to  denote  the
geodesic distance $\arccos x\cdot y$ between $x$ and $y$.
Moreover, we denote by $ B(x,r):=\{ y\in\sph: \ d(x,y)\leq r\}$
the spherical cap with center $x\in\sph$ and radius $r\in
(0,\pi)$, and  $M(f)$ the usual Hardy-Littlewood maximal function
on
 $\sph$:
$$ M(f) (x):=\sup_{0<r\leq \pi} \f 1{ |B(x,r)|} \int_{B(x,r)}
|f(y)|\,d\sa(y),\     \   x\in\sph,   \     \  f\in L(\sph).$$ We
let $\mathcal{S}\equiv\mathcal{S}(\sph)$ denote the set of
indefinitely differentiable functions on $\sph$ endowed with the
usual test function topology and let $\mathcal{S}'\equiv
\mathcal{S}'(\sph)$ be the dual of $\mathcal{S}$.  $\mathcal{S}$
is called the space of test functions and $\mathcal{S}'$ the space
of distributions.
  We denote by  $\langle f, \vi\rangle$ the  pairing between a
distribution $f\in\mathcal{S}'$ and a test function
$\vi\in\mathcal{S}$.   Throughout the paper, the notation $\# \Ld$
denotes the cardinality of a given finite set $\Ld$,  the letter
$C$ denotes a general positive constant depending only on the
parameters indicated as subscripts, and  the notation $A\asymp B$
 means that there are two inessential positive constants $C_1$,
$C_2$ such that $C_1 A\leq B \leq C_2 A$.

For each nonnegative integer $k$,
 let  $\mathcal{H}_k$ denote  the space of spherical harmonics
 of degree $k$  on $\sph$  and  $Y_k$   the orthogonal
 projection of $L^2(\sph)$ onto $\mathcal{H}_k$.  As is well
 known, for $f\in L^2(\sph)$,
\begin{equation}\label{1-1}
Y_k(f) (x) = \int_{\sph} f(y) P_k (x\cdot y)\, d\sa(y),\    \
x\in\sph,
\end{equation} where for $x=(x_1,x_2,\ldots, x_d)$,
$y=(y_1,y_2,\ldots, y_d)\in\mathbb{R}^d$, $x\cdot y =x_1y_1+x_2y_2
+\ldots x_dy_d$,
\begin{equation}\label{1-2} P_k(t):=\f{2k+d-2}{d-2} P_k^{\f{d-2}2}
(t),\    \  t\in[-1,1]
\end{equation}
and  $P_k^{\f{d-2}2}  (t)$ denotes   the usual ultraspherical
polynomial of order $\f{d-2}2$ normalized by
 $ P_k^{\f{d-2}2} (1)=\df{\Ga(k+d-2)}{\Ga(d-2) \Ga( k+1)}.$  ( For
  precise definition of  ultraspherical polynomials, we
 refer to [Sz,p. 81].)
 Evidently, the formula ($\ref{1-1}$) allows  us to extend the definition of $Y_k(f)$ to contain  all
  distributions $f\in\mathcal{S}'(\sph)$.  For an integer $n\ge 0$ we denote
 by $ \Pi_n$ the space of all spherical polynomials of degree at
 most $n$ on $\sph$  (i.e.,  polynomials in $d$-variables of
 total  degree at most $n$ restricted to $\sph$).    It is well
 known that $\displaystyle \Pi_n= \bigoplus_{k=0}^{n} \mathcal {H}_k$
and $\text{dim}\  \Pi_n =\dsum_{k=0}^n\text{dim }\
\mathcal{H}_k\asymp n^{d-1}.$  For more information on
 spherical harmonics, we refer to [SW, Ch IV].

\subsection{Construction of polynomial frames on $\sph$.}
Various nonstationary  wavelets or  frames have been constructed
on the sphere by many authors ( see [ADS, AV, DDSW, FGS, G, MNPM,
NW]). Although each of these wavelets or frames  has its own
advantages, to the best of our knowledge none of them has been
shown  useful in the  characterizations  of  classic function
spaces on the sphere. In this subsection, we will construct a
polynomial frame on $\sph$, for which  many function spaces on the
sphere $\sph$, such as $L^p$, $H^p$ and Besov spaces, can be
characterized in terms of the coefficients in the  wavelet
decompositions.

 Our construction is motivated by the approach
taken in the pioneer work [MNPM], where the authors used positive
cubature formulae to introduce a class of polynomial frames
suitable for analyzing data on the sphere.  It
 is   based on the following
 theorem:

 \vspace{5mm}
 {\bf Theorem A.}\    \  {\it There exists a constant $\ga>0$
 depending only on $d$ such that for  any
 integer $N>0$ and any finite set $ \{\xi_k\}_{k\in\Og} $ of
 distinct points $\xi_k\in\sph$ satisfying
 $$ \min_{\sub{i,j\in\Og\\ i\neq j}} d(\xi_i, \xi_j) \ge {\ga}/N,\
    \  \text{and}\    \  \max_{x\in\sph} \min_{j\in\Og} d(x, \xi_j) <{\ga}/N,$$
there exists a set of  numbers $ 0\leq a_{N,k} \leq C_d
N^{-(d-1)},\ \ k\in\Og$ such that for any  $f\in  \Pi_N$, $0<p\leq
\infty$ and  $t\ge 0$,
\begin{align} \int_{\sph}& f(y)\, d\sa(y) =\sum_{ k\in\Og}
a_{N,k} f(\xi_k),\label{1-3}\\
 \|f\|_p& \asymp \begin{cases} \( \f1{N^{d-1}}\dsum_{k\in\Og}
(N^{d-1}a_{N,k})^t  |f( \xi_k)|^p\)^{\f1p},&\   \  \text{if $0<p<\infty$,}\\
\dmax_{k\in\Og} \  \Bl[(N^{d-1}a_{N,k})^t|f(\xi_k)|\Br], &\     \
\text{ if $ p=\infty$,}
\end{cases}\label{1-4}
\end{align} where the constants of equivalence  depend only
on $d$ and $p$ when $p$ is small, and in sequel, we employ the
slight abuse of notation that $0^0=1$. }

Theorem A under the restriction $t\leq \min\{p,1\}$  was obtained
in [BD, Theorem 3.1].  It will be shown in Section 5 that the
restriction $t\leq \min\{p,1\}$ is in fact not necessary.

 An equality like ($\ref{1-3}$) with nonnegative coefficients is called  a positive  cubature formula
 while an equivalence   like ($\ref{1-4}$) is called an  Marcinkiewicz-Zygmund
 (MZ) type
inequality.    It should be pointed out that positive cubature
formulae and  MZ inequalities for $1\leq p\leq \infty$    based on
function values at scattered sites on $\sph$ were first
established in the fundamental paper  [MNW].

Now to each  integer $j>0$ we assign a finite set  $\{ x_{j,k}:\ \
k\in\Ld_j^d\}$  of distinct points $x_{j,k}\in\sph$ satisfying
$$ \min_{\sub{k,k'\in\Ld_j^d\\  k\neq k'}}  d( x_{j,k}, x_{j, k'})\ge \f \ga { 2^{j+4}},\     \    \text{and}\     \
\max_{x\in\sph} \min_{k\in\Ld_{j}^d} d( x,  x_{j,k}) <\f{\ga} {
2^{j+4}},$$ with $\ga$ as in Theorem A.   Evidently, $\#\Ld_j^d
\asymp 2^{j(d-1)}$ and by Theorem A,  there exists a set of
numbers $0\leq \ld_{j,k}\leq C_d  2^{-j(d-1)},\   \ k\in\Ld_j^d$
such that for any $0<p\leq \infty$,  $ 0\leq t \leq \min\{p,1\}$
and any $f\in\Pi_{ 2^{j+4}}$,
\begin{equation}\label{1-5}
 \int_{\sph} f(y)\,
d\sa(y)=\dsum_{ k\in\Ld_j^d} \ld_{j,k} f(x_{j,k}),\end{equation}
\begin{equation}\label{1-6} \|f\|_p\asymp
\begin{cases} \( \df1{2^{j(d-1)}}\dsum_{ k\in\Ld_j^d}
(2^{j(d-1)}\ld_{j,k})^t  | f(x_{j,k})|^p\)^{\f1p},\ \ &\text{ if
$0<p<\infty$,}\\
\dmax_{ k\in\Ld_j^d} \Bl[(2^{j(d-1)}\ld_{j,k})^t
|f(x_{j,k})|\Br],&\ \ \text{if  $
p=\infty$,}\end{cases}\end{equation} where the constants of
equivalence  depend only on $d$ and $p$ when $p$ is small. For
convenience,  we also  set $\Ld_0^d=\{0\}$, $\ld_{0,0}=1$ and take
$x_{0,0}$ to be any fixed point on $\sph$.

Let $\phi$ be a nonnegative $C^\infty$- function  on $\mathbb{R}$
supported in $\{x\in\mathbb{R}:\    \f 12\leq |x|\leq 2\}$ and
satisfying \begin{equation}\label{1-7}
 \dsum_{j=-\infty} ^\infty \( \phi(
2^{-j} x)\)^2 =1,\     \  \text{for all $x\neq 0$}.\end{equation}
Together with $\phi$ we define  a sequence of functions
$$  \psi_{j,k}(x):=\sqrt{\ld_{j,k}} G_j( x\cdot
x_{j,k}),\    j\ge 0,\    k\in\Ld_j^d, $$ where
\begin{equation}\label{1-8} G_0(t)=1,\ \ G_j( t)=\sum_{ k=[2^{j-2}]}^{
2^{j}} \phi( \f k { 2^{j-1}}) P_k(t),\
    \     t\in[-1,1],\   \   j\ge 1\end{equation}
    and $P_k(t)$ is  defined by ($\ref{1-2}$).
    We index these functions by the spherical caps
    $B(x_{j,k}, 2^{-j}\pi)$.  Thus with $B=B(x_{j,k}, 2^{-j}\pi)$ we let
    \begin{equation}\label{1-9} \psi_B(x)=\psi_{j,k}(x) =\sqrt{\ld_{j,k}} G_j( x\cdot
    x_{j,k}).\end{equation}
    We also use the notation $\mathcal{B}_j$ to denote the set of
    spherical caps $B(x_{j,k}, 2^{-j}\pi)$, ( $ k\in\Ld_j^d$) and
    $\mathcal{B}$ to denote the union of the $\mathcal{B}_j$,
    $j\ge 0$.  Thus, associated with  each distribution $f$ on
    $\sph$, there is a series $
    \dsum_{B\in\mathcal{B}} \lb f,\psi_B\rb \psi_B,$ and   moreover, for
    each spherical polynomial $f$, we have
        \begin{equation*}\label{1-1n} f(x)=  \sum_{B\in\mathcal{B}} \lb f,\psi_B\rb
        \psi_B(x),\end{equation*}
    with only a finite number of nonzero coefficients $\lb
    f,\psi_B\rb$,  as can be easily verified.  We will keep the
    above notations for the rest of the paper.

We have three   purposes  in   this  paper. First, we want to
study the unconditional  convergence of the series
$\dsum_{B\in\mathcal{B}} \lb f,\psi_B\rb \psi_B(x)$. Second, we
want to characterize classic function spaces on $\sph$, such as
$L^p$, $H^p$ and Besov spaces, using  the coefficients $\lb
f,\psi_B\rb$, ($B\in\mathcal{B}$). Third, we wish to investigate a
related nonlinear $m$-term approximation problem.   Our main
results will be summarized in the next subsection.

\subsection{Summary of main results}
To state our results in a unified manner, we  identify
$H^\infty(\sph)$ with $C(\sph)$, and note  that (see [C]) for $
1<p< \infty$,  the Hardy space $H^p(\sph)$ coincides with the
Lebesgue space $L^p(\sph)$ and $\|f\|_p\asymp \|f\|_{H^p}$ with
the constants of equivalence depending only on $d$ and $p$.

Our first result gives a characterization of the spaces $H^p$,
$0<p<\infty$.

\begin{thm} If $0<p<\infty$ and $f\in H^p(\sph)$, then
$\dsum_{B\in\mathcal{B}} \lb f,\psi_B\rb \psi_B $ converges
unconditionally  to  $f$ in the $H^p$-metric, and moreover
$$\|f\|_{H^p}\asymp  \Bl\| \(
\sum_{ B\in\mathcal{B}} |\langle f,\psi_B\rangle|^2 |\psi_B |^2
\)^{\f12} \Br\|_p\asymp  \Bl\|\( \sum_{ B\in\mathcal{B}} |\langle
f, \psi_B\rangle |^2 |B|^{-1} \chi_B\)^{\f12}\Br\|_p,$$ with the
constants of equivalence depending only on $d$, $p$ and $\phi$. In
addition, if $0<p<\infty$ and $\{a_B\}_{B\in\mathcal{B}}$ is a
sequence of complex numbers such that either
$\(\dsum_{B\in\mathcal{B}} |a_B|^2|\psi_B|^2\)^{\f12}\in
L^p(\sph)$ or $\(\dsum_{B\in\mathcal{B}} |a_B|^2|B|^{-1}
\chi_B\)^{\f12}\in L^p(\sph)$, then  $\dsum_{B\in\mathcal{B}} a_B
\psi_B$ converges unconditionally  to some  distribution $f$ in
the $H^p$-metric, and moreover
$$\|f\|_{H^p} \leq C_1 \Bl\|\(\sum_{B\in\mathcal{B}}
|a_B|^2|\psi_B|^2\)^{\f12}\Br\| _p \leq C_2 \Bl\|
\(\sum_{B\in\mathcal{B}} |a_B|^2|B|^{-1} \chi_B\)^{\f12}\Br\|
_p,$$ with the constants $C_1$ and $C_2$  depending only on $d$,
$p$ and $\phi$.
\end{thm}

Our next result concerns  a  characterization  of the  Besov
spaces.  For $\al>0$ and $ 0<p,\tau\leq \infty$, we define
$$|f|_{B_\tau^\al (H^p)} :=\( \sum_{k=0}^\infty (k+1)^{\al\tau-1}
(E_k(f)_{_{H^p}})^\tau\)^{\f1\tau}+|\lb f, 1\rb|,$$
 with the
usual change when $\tau=\infty$, where
$$ E_n(f)_{_{H^p}}:=\inf\Bl\{ \|f-g\|_{H^p}:\     \
g\in\Pi_n\Br\},\   \  n\in\mathbb{Z}_+,$$ and  define the Besov
space $B_\tau^\al(H^p)$ to be a linear space  of distributions on
$\sph$    endowed with the quasi-norm $|\cdot|_{B_\tau^\al(H^p)}$.
 We point out that  Besov spaces on the sphere
were  introduced and investigated by Nikol'skii,  Lizorkin and
Rustamov  in a series of papers  ( see [R] and the references
there).  These spaces can be equivalently  characterized using the
K-functionals or moduli of smoothness on the sphere ( see [R]).

\begin{thm}
For $\al>0$, $ 0<p,\tau\leq \infty$ and $f\in B_\tau^\al (H^p)$,
we have
\begin{equation}\label{1-10} |f|_{B_\tau^\al (H^p)} \asymp  \(
\sum_{j=0}^\infty 2^{-j(d-1)\tau(\f 1p-\f12-\f\al{d-1})} \(
\sum_{k\in\Ld_j^d} |\langle f,\psi_{j,k}\rangle|^p \)^{\f \tau p}
\)^{\f1\tau},\end{equation} with the usual change when $p=\infty$
or $\tau=\infty$, where the constants of equivalence depend only
on $\al,\tau,p$ and $\phi$.  In addition, if $\{a_{j,k}:\    \
j=1,2,\ldots, \   k\in\Ld_j^d\}$ is a sequence of complex numbers
such that
$$ \( \sum_{j=0}^\infty 2^{-j(d-1)\tau(\f 1p-\f12-\f\al{d-1})}
\( \sum_{k\in\Ld_j^d} | a_{j,k}|^p \)^{\f \tau p} \)^{\f1\tau}
<\infty,$$ then  the series $\dsum_{j=0}^\infty\dsum_{k\in\Ld_j^d}
a_{j,k} \psi_{j,k}$ converges unconditionally  to some $f\in
B_\tau^\al (H^p)$ in the $H^p$-metric and moreover,
\begin{equation} \label{1-11}|f|_{B_\tau^\al (H^p)} \leq C_{p,\al,\tau,
\phi}\( \sum_{j=0}^\infty 2^{-j(d-1)\tau(\f 1p-\f12-\f\al{d-1})}
\( \sum_{k\in\Ld_j^d} | a_{j,k}|^p \)^{\f \tau p}
\)^{\f1\tau},\end{equation} with the usual change when $p=\infty$
or $\tau=\infty$.\end{thm} Of particular interest are the spaces
$B_\tau^\al(H^\tau)$, for which we have
\begin{equation}\label{1-12}|f|_{B_\tau^\al (H^\tau)} \asymp  \(\sum_{B\in\mathcal{B}} |\langle
f,\psi_B\rangle|^\tau |B|^{1-\f \tau2-\f{\al \tau}{d-1}}    \)^{\f
1\tau},\end{equation} with the usual change when $\tau=\infty$, on
account of Theorem 1.2.

We point out that for the usual Euclidean space $\mathbb{R}^d$,
results similar to Theorems 1.1 and 1.2 can be found in [FJ],[HW,
Chapter 7] and [DJP].

Finally,  we state our results on nonlinear approximation.  For
$0<p\leq \infty$,  $f\in H^p$ and an integer $n>0$, we denote by
$\Ga_n\equiv \Ga_{n,f,p}$ a set of $n$ spherical caps
$B\in\mathcal{B}$ such that $$\min_{B\in\Ga_n} |\langle
f,\psi_B\rangle||B|^{\f 1p-\f12} \ge \max_{B\in \mathcal{B}
\setminus \Ga_n}|\langle f,\psi_B\rangle||B|^{\f 1p-\f12},$$ and
define the greedy type algorithm $G_n^p(f)$ by
$$ G^p_n(f) =
\sum_{B\in\Ga_n} \langle f,\psi_{B}\rangle \psi_{B}.$$ Such an
algorithm is well defined, as was shown in  [T1, remark 1.1]. We
refer to the impressive survey paper [T2] for the background
information of greedy algorithm.

\begin{thm}
For $\al>0, \   0<p\leq \infty$, $\tau=\( \f\al{d-1}+\f1p\)^{-1}$
, $f\in H^p$ and an integer $n>0$,
\begin{equation}\label{1-13} \|f-G^p_n(f)\|_{H^p} \leq C_{p,\al,\phi} n^{-\f\al{d-1}}
|f|_{B_\tau^\al(H^\tau)}.\end{equation}
\end{thm}

\begin{thm}
For $\al>0$, $0<p\leq \infty$, $\tau=(\f 1p+\f\al{d-1})^{-1}$,
$f\in H^p$ and an integer $n>0$,
\begin{equation}\label{1-14}\|G^p_n(f)\|_{B_\tau^\al(H^\tau)} \leq C_{p,\al,\phi} n^{\f\al {d-1}}
\|f\|_{H^p}.\end{equation}   In addition, if $f\in H^p$ has a
representation $f(x)=\dsum_{B\in\mathcal{B}} a_B\psi_B(x)$ with at
most $n$ nonzero coefficients  $a_B$, then
\begin{equation}\label{1-15}|f|_{B_\tau^\al (H^\tau)} \leq C_{p,\al,\phi} n^{\f\al{d-1}}
\Bl\| \max_{B\in \mathcal{B}}\Bl[ |a_B| |B|^{-\f12} \chi_B(\cdot)
\Br]\Br\|_p.\end{equation} \end{thm}

 The inequality ($\ref{1-13}$) is a direct theorem of
 approximation (the Jackson inequality) while the inequality
 ($\ref{1-14}$) is an inverse theorem ( the Bernstein inequality).
  Once ($\ref{1-13}$) and ($\ref{1-14}$) are established, then by the
  standard method (see [DP]), the following characterization
  result holds for $0<\be<\al$:
  $$ \sum_{n=1}^\infty \Bl[ 2^{n\be/(d-1)} \|f-G_{2^n} ^p
  (f)\|_{H^p}\Br]^\tau <\infty  \     \   \Longleftrightarrow   \
    \   f\in B_\tau^\be (H^\tau),$$
    where  $f\in H^p$, $0<p\leq \infty$ and $\tau= ( \be /(d-1)
    +1/p)^{-1}$.
    For results on the  nonlinear approximation associated with   wavelet
    decomposition in  $L^p(\mathbb{R}^d)$  ( $0<p\leq \infty$),
we refer to [DJP], [DPY] and [Jia].

We organize this paper as follows. Section 2 contains three lemmas
which will be used frequently in the proofs of our main results.
We prove Theorem 1.1  in Section 3 for the case $1<p<\infty$, and
in Section 4 for the case $0<p\leq 1$.  The proof of Theorem 1.2
and those of Theorems 1.3 and 1.4 are given in Sections 5 and 6
respectively.

Finally, we point out that  our paper does not give any effective
algorithms for spherical wavelets. For results in this direction,
we refer the reader to the  papers  [SS1], [SS2] and [KL]. We also
note that characterizations of function spaces using wavelets on
stratified Lie groups had been considered in [L]. (The author
would like to thank an  anonymous referee very much for kindly
pointing out these  references ([SS1], [SS2], [KL], [L]) to him.)
For more recent work on spherical frames, we refer to the
impressive  paper [NPW] by F. J. Narcowich, P. Petrushev  and J.
D. Ward,  and also the nice  survey paper [MP] by   H. N. Mhaskar
and  J. Prestin.

\section{Three useful  lemmas}

This section contains three lemmas that will be useful in the
proofs of the  main results in this paper. The first two lemmas (
Lemmas 2.1 and 2.2) are in essence known, while the last one
(Lemma 2.3) is new and will be of fundamental importance.

For the statement of Lemma 2.1,  we define, for
$f\in\mathcal{S}'(\sph)$,
        \begin{equation}\label{2-1} \sa_j(f)(x) := \lb f, G_j(x\cdot)\rb,\      \
       x\in \sph,\      \   j=0,1,2,\ldots,\end{equation}
       where $G_j$ is defined by ($\ref{1-8}$).
      \begin{lem}
 For $1<p<\infty$ and $f\in L^p(\sph)$,
              \begin{equation}\label{2-2}
       \Bl\|\(\sum_{j=0}^\infty
       |\sa_j(f)|^2\)^{\f12}\Br\|_p\leq C_p \|f\|_p.
       \end{equation}
In addition, if $f$ is a spherical polynomial, then
\begin{equation}\label{2-3}\sum_{j=0}^\infty \sa_j \circ \sa_j (f)(x) =f(x),\     \
\forall\ x\in\sph,\end{equation} with only a finite number of
nonzero terms.
\end{lem}

\begin{proof}
The inequality ($\ref{2-2}$) is a simple consequence of the
well-known H\"ormander-type multiplier theorem for spherical
harmonics (see [S, Theorem 2]), while the identity ($\ref{2-3}$)
follows directly from ($\ref{1-7}$) and the definition.\end{proof}

To state our next lemma, we suppose  $\vi$ is  a $C^\infty$
-function on $[0,\infty)$ supported in $[0,2]$ and  equal to a
constant on $[0,\f12]$, and  define
\begin{equation}\label{2-4} K_{N,\vi}(t):=\sum_{k=0}^\infty \vi(\f kN) P_k(t),\
\    t\in[-1,1],\    \  N=1,2,\ldots, \end{equation} with $P_k(t)$
as defined in ($\ref{1-2}$).  Then, with these notations, we have
\begin{lem}
For   $\ta\in [0,\pi]$  and any  positive integer $\ell$,
$$ |K_{N,\vi}^{(i)}(\cos\ta)| \leq C_{\vi, \ell,i} N^{d-1+2i} \min \Bl\{ 1,
(N\ta)^{-\ell}\Br\},\   \   i=0,1,\ldots, \  N=1,2,\ldots,$$ where
$ K_{N,\vi}^{(0)}(t)=K_{N,\vi}(t)$, $K_{N,\vi}^{(i)}(t)=\(\f
d{dt}\)^i\{ K_{N,\vi}(t)\},\ \
 i\ge 1.$
\end{lem}

We note that for the usual Ces\`aro kernel
 $$\sa_N^\da(t):=\sum_{k=0}^N \f{ \Ga(N-k+\da+1) \Ga(N+1)}{\Ga(N-k+1)\Ga(N+\da+1)} P_k(t)$$  of
order $\da>d-1$  it is well-known (see [BC]) that
$$\Bl|\sa_N^\da (\cos\ta) \Br|\leq C_\da  N^{d-1} \min \Bl\{ 1,
(N\ta)^{-d}\Br\} $$ and  the order $(N\ta)^{-d}$ on the right-hand
side of  this last inequality cannot be further improved.  The
significance point of  Lemma 2.2 is that the positive integer
$\ell$ can be chosen  as big  as we like, which  will play a very
important role when we deal with the case  $0<p<1$ in the later
sections.

The proof of Lemma 2.2 is contained in [BD, Lemma 3.3].

The following lemma will be of fundamental importance in the
proofs of our main results.

\begin{lem}
Suppose $f$ is a spherical polynomial of degree at most $N$,
(i.e.,  $f\in\Pi_N$ ).  Then for any  $\be>0$ and $x\in\sph$,
\begin{equation} \label{2-5}f_{\be, N}^\ast (x) \leq C_\b  \( M ( |f|^{\f
1\be})(x) \)^\be,\end{equation} where
\begin{equation}\label{2-6} f_{\be, N} ^{\ast} (x) :=\sup_{ y\in\sph} \f{
|f(y)|} { \( 1+ N d(x,  y)\)^{\be(d-1)} }.\end{equation}
\end{lem}
\begin{proof}
Let $\eta$ be a $C^\infty$ -function on $[0,\infty)$  supported in
$[0,2]$ and  equal to $1$ on $[0,1]$ and let  $K_{N,\eta}$ be
defined  by ($\ref{2-4}$) with $\vi$ replaced by $\eta$.   Then,
clearly, for $f\in \Pi_N$, \begin{equation}\label{2-7}  f(y)-f(z)
= \int_{\sph} f(u) \(K_{N,\eta}(y\cdot u)-K_{N,\eta}(z\cdot u)\)\,
d\sa(u),\ \ y,z\in\sph.\end{equation} For simplicity, we set, for
$y,u\in\sph$ and $\da\in (0,\f14)$,
$$ A_{N,\da} (y,u) := \max_{z\in B(y,\f\da N)} | K_{N,\eta}
(y\cdot u) - K_{N,\eta} (z\cdot u)|.$$ Then, using Lemma 2.2 with
$i=0,1$, it's easy to verify that for any integer $\ell>0$,
\begin{equation}\label{2-8}
A_{N,\da} (y,u)\leq C_{\ell,\eta}
\begin{cases}
N^{d-1},&\   \text{if $\ta\in [0,\f{4\da}N]$,}\\
\da N^{d-1} \min\{ 1, (N\ta)^{-\ell}\},&\    \  \text{if $\ta\in
[\f{4\da}N,\pi]$,}
\end{cases}\end{equation}
where $\ta=d(y, u)$. Thus, using ($\ref{2-7}$) and ($\ref{2-8}$)
with $\ell= [(\b+1)(d-1)]+2$,  we obtain that  for $x,y\in\sph$
and $\da\in (0,\f14)$,
\begin{align*} &\sup_{z\in B(y,\f\da N)}\f{|f(y)-f(z)|}{ (
1+ N d(x, y))^{\be(d-1)} }\\
&\leq  f_{\be, N}^\ast(x) \int_{\sph} \( \f{ 1+ N d( x, u)}{ 1+ N
d( x,  y )}\)^{\be(d-1)} A_{N,\da} (y,u)\,
d\sa(u) \\
&\leq 2^{\b(d-1)}   f_{\be, N}^\ast(x) \Bl[ \int_{\sph} A_{N,\da}
(y,u) \, d\sa(u)  +  \int_{\sph}  ( 1+ N d(
y,u))^{\be(d-1)}A_{N,\da} (y,u)\,
d\sa(u)\Br]\\
&\leq C_{\eta,\b} \da  f_{\be, N}^\ast(x) .\end{align*}   This
implies that for any $x, y\in\sph$ and $\da\in (0,\f14)$,
\begin{align*}
&\f{|f(y)|^{\f 1\be} } { ( 1+ N d( x, y))^{d-1} } \leq
\f{2^{\f1\b} \displaystyle\inf_{ z\in B( y,\f \da N)} |f(z)|^{\f
1\be} } { ( 1+ N d(x,  y))^{d-1} }+ \(2C_{\eta,\b}\da f_{\be, N}^\ast(x) \)^{\f 1\be} \\
\leq & \f{C_\b \(\f N\da\)^{d-1} } { (1+ N d( x,  y))^{d-1}
}\int_{ B(y,\f \da N )} |f(z)|^{\f 1\be} \, d\sa(z) +
\(2C_{\eta,\b}\da f_{\be,
N}^\ast(x) \)^{\f 1\be} \\
\leq& C_\b   \da^{-(d-1)} M( |f|^{\f 1\be})(x) + \(2C_{\eta,\b}\da
f_{\be, N}^\ast(x) \)^{\f 1\be}.\end{align*} So, taking the
supremum over all $y\in\sph$, we deduce
$$ f_{\be,N}^\ast (x) \leq C_\b'  \da^{ -(d-1)\be } \( M ( |f|^{\f
1\be})(x)\)^\be + 2^{\b+1} C_{\eta,\b}\da f_{\be, N}^\ast(x),$$
 which implies ($\ref{2-5}$) by taking
 $\da=(2^{\b+2}C_{\eta,\b})^{-1}$.
This completes the proof.
\end{proof}

\section{Proof of Theorem 1.1 for $1<p<\infty$}

  We need the
following  lemma, which means that for each $j\in\mathbb{Z}_+$ and
$k\in\Ld_j^d$  the  function $\psi_{j,k}$ is highly localized in
the spherical cap $B(x_{j,k}, 2^{-j-1})$.

\begin{lem}
For any $r>0$,  $x\in\sph$,  $j\in\mathbb{Z}_+$ and $k\in\Ld_j^d$,
\begin{equation} \label{3-1}  g_{j,k}(x)\leq C_{d,\phi}  |\psi_{j,k}(x)|\leq C_{d,\phi,r}  \( M(|g_{j,k}|^r)(x)\)^{\f 1r},
\end{equation} where $g_{j,k}(x) =2^{j(d-1)}
\sqrt{\ld_{j,k}}\chi_{_{B(x_{j,k}, 2^{-j-1})}}(x)$.
\end{lem}

\begin{proof} When $j=0$, Lemma 3.1 is trivial. So we may assume
$j>0$. Recall that for $j>0$,
$$ \psi_{j,k}(x)=\sqrt{\ld_{j,k}} G_j(x\cdot x_{j,k})=\sqrt{\ld_{j,k}} \sum_{ v=0}^{2^j} \phi( \f { v}{ 2^{j-1}} ) P_v( x\cdot x_{j,k})
 ,$$
 where $\phi$ is a nonnegative $C^\infty$-function supported in $\{
 x\in\mathbb{R}:\   \f12\leq |x|\leq 2\}$ and satisfying ($\ref{1-7}$), and
 $P_v$ is defined by ($\ref{1-2}$). Since (see [Sz, (7.32.2) ]) $
 \dmax_{t\in[0,\pi]} |P_v(\cos t)| =P_v(1)\asymp v^{d-2}$, it
 follows that
 $  \dmax_{t\in [0,\pi]} | G_j(\cos t)|
 = G_j(1) \asymp 2^{j(d-1)}.$
 Hence, using Bernstein's inequality for trigonometric
 polynomials, we obtain that for $t\in [0,\f 1{ 2^{j+1}}]$,
 $$\sqrt{\ld_{j,k}}   G_j(\cos t)\geq \f12 \sqrt{\ld_{j,k}}
 G_j(1)\asymp 2^{j(d-1)} \sqrt{\ld_{j,k}},$$
 and the inequality
 $g_{j,k}(x)\leq C_{d,\phi}  |\psi_{j,k}(x)|$ then  follows.

For the proof of the inequality
\begin{equation}\label{3-2}
|\psi_{j,k}(x)|\leq C_{d,\phi,r}  \( M(|g_{j,k}|^r)(x)\)^{\f
1r},\end{equation} we use Lemma 2.2 to obtain that
 $$ |\psi_{j,k}(x)| \leq C_{d,\phi,r}  \sqrt{\ld_{j,k}} 2^{j(d-1)} \min\{
 1, (2^j\ta)^{-\f{d-1}r}\},$$
 where $\ta=d(x,  x_{j,k})$. Therefore, if $ 0\leq \ta=d( x,  x_{j,k})\leq 2^{-(j+1)}$
 then since
 $B(x_{j,k},2^{-j-1})\subset B(x, 2^{-j})$,
 \begin{align*}
 |\psi_{j,k}(x)|
  \leq & C_{d,\phi,r}
\sqrt{\ld_{j,k}} 2^{j(d-1)} \( \f 1{ | B(x, 2^{-j})|} \int_{ B(x,
2^{-j})} \chi_{_{B(x_{j,k}, 2^{-j-1})}} (y)\,
d\sa(y) \)^{\f 1r}\\
 \leq &  C_{d,\phi,r} \( M(|g_{j,k}|^r)(x)\)^{\f
1r};   \end{align*} and  if $\ta= d(x, x_{j,k})
>2^{-j-1}$ then since $ B(x_{j,k}, 2^{-j-1})\subset B(x, 2\ta)$,
\begin{align*}
|\psi_{j,k}(x)|&\leq C_{d,\phi,r} \sqrt{\ld_{j,k}} 2^{j(d-1)} (
2^j\ta)^{-(d-1)/r}\\
& \leq C_{d,\phi,r} \sqrt{\ld_{j,k}} 2^{j(d-1)} \( \f 1{ | B(x,
2\ta)|}\int_{B(x,2\ta)} \chi_{_{B(x_{j,k}, 2^{-j-1})}} (y) \,
d\sa(y) \)^{\f1r} \\
&\leq C_{d,\phi,r}  \( M(|g_{j,k}|^r)(x)\)^{\f 1r}.
\end{align*}
In either case, we have the desired estimate ($\ref{3-2}$), and
the proof is therefore  complete.
 \end{proof}

{\it Proof of Theorem 1.1 for $1<p<\infty$.}\    \   First, we
show that for any sequence  $\{a_B\}_{B\in\mathcal{B}}$ of complex
numbers,
 \begin{equation}\label{3-3} \Bl\| \(
\sum_{ B\in\mathcal{B}} |a_B|^2 |\psi_B|^2 \)^{\f12} \Br\|_p\leq
C_{p, d,\phi}  \Bl\|\( \sum_{ B\in\mathcal{B}}  |a_B|^2 |B|^{-1}
\chi_B\)^{\f12}\Br\|_p.\end{equation} Indeed, since  $0\leq
\ld_{j,k}\leq C_d 2^{-j(d-1)}$,  by Lemma 3.1 applied to $r=1$ it
follows that for any $B\in\mathcal{B}$, $ |\psi_B| \leq C_{d,\phi}
|B|^{-\f12} M(\chi_{B}),$  which together with the well-known
Fefferman-Stein inequality implies the inequality ($\ref{3-3}$).

Second, we show that for any sequence  $\{a_B\}_{B\in\mathcal{B}}$
of complex numbers and any finite subset $\mathcal{F}\subset
\mathcal{B}$,
\begin{equation}\label{3-4}
\|\sum_{B\in\mathcal{F}} a_B\psi_B\|_p \leq C_{p,d,\phi}
\Bl\|\(\sum_{B\in\mathcal{F}} |a_B|^2 |\psi_B|^2 \)^{\f12}\Br\|_p.
\end{equation}
Once ($\ref{3-4}$) is proved, then by a standard argument, we
deduce that the series $\dsum_{B\in\mathcal{B}} a_B\psi_B$
converges unconditionally in the space $L^p$ provided that
$\(\dsum_{B\in\mathcal{B}} |a_B|^2 |\psi_B|^2 \)^{\f12}\in L^p$.
This together with ($\ref{3-3}$) will imply the second assertion
of Theorem 1.1  in the case $1<p<\infty$.

For the proof of ($\ref{3-4}$), we define
\begin{equation}\label{3-5} \sa_j^{\ast\ast} (f)(x) =\sup_{y\in\sph} \f{ |\sa_j(f)(y)|}
{ ( 1+2^j d(x, y))^{d-1}},\end{equation} where $\sa_j(f)$ is
defined by ($\ref{2-1}$). For $j\in\mathbb{Z}_+$ and
$k\in\Ld_j^d$, we write  $B_{j,k}=B( x_{j,k}, 2^{-j}\pi)$,
$B'_{j,k}=B( x_{j,k}, 2^{-j-1})$, and   we define
$a'_{j,k}=a_{_{B_{j,k}}}$ if $B_{j,k} \in\mathcal{F}$, and $=0$
otherwise. Also, we set \mbox{$h=\dsum_{B\in\mathcal{F}} a_B
\psi_B$}. Let $g\in L^{p'}$ be such that $\|g\|_{p'}=1$ and
$\|h\|_p=\dint_{\sph} f(x) g(x)\, d\sa(x)$, where $p'=\f{p}{p-1}$.
We observe that
$$
\Bl|\lb g,\psi_{j,k}\rb\Br|=\sqrt{\ld_{j,k}} |\sa_j(g)(x_{j,k})|
\leq C_d  2^{j(d-1)} \sqrt{\ld_{j,k}} \int_{B_{j,k}'}
\sa_j^{\ast\ast}(g)(x)\, d\sa(x).$$ Therefore
\begin{align*}
\|h\|_p&= \sum_{j=0}^\infty\sum_{k\in\Ld_j^d}  a'_{j,k}  \langle
\psi_{j,k},g\rangle\leq C_d \sum_{j,k}  |a'_{j,k}|
\sqrt{\ld_{j,k}} 2^{j(d-1)} \int_{B_{j,k}'}
\sa_j^{\ast\ast}(g) (x) \, dx\\
&\leq  C_d \Bl\|\(\sum_{j,k} |a'_{j,k}|^2 \ld_{j,k} 2^{2j(d-1)}
\chi_{ _{B_{j,k}'}} \)^{\f12}\Br \|_p
\Bl\|\(\sum_{j=0}^\infty\sum_{k\in\Ld_j^d}
|\sa_j^{\ast\ast}(g)|^2\chi_{_{B_{j,k}'}}\)^{\f12} \Br\|_{p'}\\
&\leq C_d \Bl\|\(\sum_{j,k} |a'_{j,k}|^2 \ld_{j,k} 2^{2j(d-1)}
\chi_{ _{B_{j,k}'}} \)^{\f12}\Br \|_p
\Bl\|\(\sum_{j=0}^\infty|\sa_j^{\ast\ast}(g)|^2\)^{\f12}
\Br\|_{p'}.
\end{align*} Invoking Lemma 2.1, Lemma 2.3 with $\b=1$
and $N=2^j$, and the Fefferman-Stein inequality, we deduce
$$\Bl\|\(\sum_{j=0}^\infty|\sa_j^{\ast\ast}(g)|^2\)^{\f12}
\Br\|_{p'}\leq C_p \Bl\|\(\sum_{j=0}^\infty|\sa_j (g)|^2\)^{\f12}
\Br\|_{p'}\leq C_p,$$ while using the first inequality in Lemma
3.1, we have
$$ \(\sum_{j,k} |a'_{j,k}|^2 \ld_{j,k} 2^{2j(d-1)}
\chi_{ _{B_{j,k}'}} \)^{\f12}\leq C_{d,\phi} \(\sum
_{B\in\mathcal{F}} |a_B|^2 |\psi_B|^2 \)^{\f12}.$$ The desired
inequality ($\ref{3-4}$) then follows.

Finally, we show the first assertion of Theorem 1.1.  We claim
that it will suffice to prove that for $f\in L^p$,
\begin{equation}\label{3-6} \Bl\|\( \sum_{ B\in\mathcal{B}}
|\langle f, \psi_B\rangle |^2 |B|^{-1} \chi_B\)^{\f12}\Br\|_p \leq
C_{p, d,\phi}  \|f\|_p.\end{equation} In fact, once ($\ref{3-6}$)
is proved, then by the second assertion of Theorem 1.1 we just
proved, it follows that  for $f\in L^p$, the series
$\dsum_{B\in\mathcal{B}} \lb f,\psi_B\rb \psi_B $ is convergent
unconditionally in $L^p$, and by the usual density argument we
must have $f=\dsum_{B\in\mathcal{B}} \lb f,\psi_B\rb \psi_B. $
This together with ($\ref{3-3}$), ($\ref{3-4}$) and ($\ref{3-6}$)
will imply the desired equivalences
$$\|f\|_{p}\asymp  \Bl\| \(
\sum_{ B\in\mathcal{B}} |\langle f,\psi_B\rangle|^2 |\psi_B |^2
\)^{\f12} \Br\|_p\asymp  \Bl\|\( \sum_{ B\in\mathcal{B}} |\langle
f, \psi_B\rangle |^2 |B|^{-1} \chi_B\)^{\f12}\Br\|_p,$$ and hence
the first assertion of Theorem 1.1.

For the proof of ($\ref{3-6}$), we  recall that $ \lb
f,\psi_{j,k}\rb =\sqrt{\ld_{j,k}} \sa_j(f) (x_{j,k})$ and $ 0\leq
\ld_{j,k} \leq C_d 2^{-j(d-1)}$.
 Thus, by the definition,  it's easy to verify that
$$
\sum_{j=0}^\infty \sum_{k\in\Ld_j^d} |\lb f,\psi_{j,k}\rb |^2
|B_{j,k}|^{-1} \chi_{_{B_{j,k}}} (x) \leq C_d \sum_{j=0}^\infty
|\sa_j^{\ast\ast} (f) (x) |^2,$$ where $\sa_j^{\ast\ast}$ is
defined by ($\ref{3-5}$). The desired inequality ($\ref{3-6}$)
then follows by Lemma 2.1, \mbox{Lemma 2.3} and the well-known
Fefferman-Stein inequality. This completes the proof.\hb

\section{Proof of Theorem 1.1 for $0<p\leq 1$}

We start with some basic definitions and facts related to the
Hardy spaces $H^p(\sph)$, $ 0<p<\infty$.   For $x\in\sph$ and
$z\in B_d:=\Bl\{ (z_1,\cdots, z_d)\in \mathbb{R}^d: \
z_1^2+\cdots+z_d^2 < 1\Br\}$,  let $p_z(x) = c_d \df {1-|z|^2} {
|z-x|^d},$ where $c_d$ is chosen so that $\dint_{\sph} p_z(x)\,
d\sa(x)=1$ for all $z\in B_d$.   $p_z$ belongs to
$\mathcal{S}(\sph)$ and is called the Poisson Kernel.   Given  a
distribution $f\in \mathcal{S}'(\sph)$, we define its radial
maximal function  by
$$P^{+}f(x) =\sup _{0\leq r <1} |\lb f, p_{_{rx}}\rb|,\    \     \   x\in\sph$$
 and its $H^p$ -quasi-norm  ( for a given $0<p<
\infty$) by $\|f\|_{H^p} =\|P^{+}f\|_p$.    For $0<p< \infty$, the
Hardy space $H^p(\sph)$  is   defined to be the space of all
distributions $f\in\mathcal{S}'(\sph)$ with $\|f\|_{H^p}<\infty$.
It is well known (see [C]) that  if $1<p< \infty$  then the Hardy
space $H^p$ coincides with the Lebesgue space $L^p$  and
$\|f\|_p\asymp \|f\|_{H^p}$ with the constants of equivalence
depending only $p$ and $d$. We will restrict ourselves to the case
$0<p\leq 1$ for the rest of this section.

For $ 0<p\leq 1$, $ 1\leq q \leq \infty$ and a nonnegative integer
$s$, a regular $(p,q,s)$ -atom centered at a point $y\in\sph$ is a
function $a$ in $L^q(\sph)$ satisfying the following three
conditions:

(i)\   \   $\text{supp}\ a\subset B(y,r)$ for some $r>0$;

(ii)\    \  $\|a\|_q\leq r^{(d-1)(\f1q-\f1p)}$;

(iii)\    \  $\dint_{\sph} a(x) p(x)\, d\sa(x)=0$ for all
$p\in\Pi_s$.\\
 An exceptional atom is a function $a$ in
$L^\infty(\sph)$ with $\|a\|_\infty\leq 1$.  Then the well-known
atomic decomposition theorem (see [C, Proposition 3.1]) states
that if $0<p\leq 1$, $ 1<q \leq \infty$, $ s\ge [ (d-1) (\f1p-1)]$
and $f\in H^p(\sph)$ then there exist a sequence $\{
c_j\}_{j=0}^\infty$ of complex numbers and a sequence
$\{a_j\}_{j=0}^\infty$ of exceptional  or regular $(p,q,s)$-atoms
such that $ \dsum_{j=0}^\infty c_j a_j$ converges to $f$ in the
space  $H^p$ and
$$
\(\sum_{j=0}^\infty |c_j|^p\)^{\f1p} \leq C_{p,d} \|f\|_{H^p}.$$

Given $ 0<p\leq 1$, a $p$-molecule centered at a point $y\in\sph$
is a function $m\in L^2(\sph)$ satisfying the following two
conditions:

(i${}'$)\    \ For some $r>0$ and $s>(d-1)(\f2p-1)$,
\begin{equation*}
\(\int_{\sph} |m(x)|^2 \( 1+ \f{d( x, y)} r\)^s \,
d\sa(x)\)^{\f12} \leq r^{-(d-1) (\f1p-\f12)};\end{equation*}

(ii${}'$)\    \
\begin{equation*}
\int_{\sph} m(x)p(x)\,d\sa(x)=0,\     \      \   \text{ for all
$p\in\Pi_{[(d-1)(\f1p-1)]}$}.\end{equation*} According to [C, p.
234], for $0<p\leq 1$, any $p$-molecule $m$ must satisfy
$\|m\|_{H^p}\leq C_{p,d}$.

For the proof of Theorem 1.1, we need the following
\begin{lem}
For $0<p\leq 1$ and $f\in H^p$, $$\Bl\| \( \sum_{ j=0}^\infty
|\sa_j(f)|^2 \)^{\f12} \Br\|_p \leq C_{d,p,\phi} \|f\|_{H^p},$$
where $\sa_j$ is defined by ($\ref{2-1}$).\end{lem}
\begin{proof} Let $0<p\leq 1$ and  $ s=[(d-1)(\f1p-1)]+2$.
Since $\(\dsum_{j=0}^\infty|\sa_j(f)|^2\)^{\f12}$ is bounded on
$L^2$,  by the atomic decomposition theorem it will suffice to
prove that  for any  regular $(p,\infty,s)$-atom $a$,
\begin{equation}\label{4-1}
\Bl\|\(\sum_{j=0}^\infty|\sa_j(a)|^2\)^{\f12}\Br\|_{p}\leq
C_{p,d,\phi}.\end{equation}

 For the proof of ($\ref{4-1}$), we suppose $a$ is a regular
$(p,\infty,s)$-atom supported in $B(x_0,r)$ for some $x_0\in\sph$
and $r\in (0,\f14)$.  Then  using  H\"older's inequality, we have
\begin{align*}
\( \int_{ B(x_0, 4r)}& \( \sum_{j=0}^\infty |\sa_j(a)(x)|^2\)^{\f
p2} \, d\sa(x)\)^{\f1p}\leq \( \int_{B(x_0,4r)} \sum_{j=0}^\infty
|\sa_j(a)(x)|^2\, d\sa(x) \)^{\f12} \(\int_{B(x_0,
4r)} d\sa(x)\)^{ \f {2-p}{2p}}\notag \\
&\leq C_{p,d,\phi}  r^{\f{(2-p)(d-1)}{2p} } \|a\|_2\leq
C_{p,d,\phi}'.\end{align*} Thus, it remains to prove
\begin{equation}\label{4-2}
\( \int_{ \sph \setminus B(x_0, 4r)}\( \sum_{j=0}^\infty
|\sa_j(a)(x)|^2\)^{\f p2} \, d\sa(x)\)^{\f1p}\leq
C_{p,d,\phi}.\end{equation} To prove this last inequality,  we
claim that for $x\in \sph\setminus B(x_0,4r)$ and  $\ell=d+2s+2$,
\begin{equation}\label{4-3}|\sa_j(a)(x)|\leq C_{p,d,\phi} 2^{j(d+2s+1)} \ta^{s+1}
r^{s-\f{d-1}p+d} \min\{ 1, (2^j\ta)^{-\ell}\},\end{equation} where
$\ta=d( x, x_0)$. Once  the claim ($\ref{4-3}$) is proved, then by
straightforward calculation,  we deduce that for $x\in
\sph\setminus B(x_0,4r)$,
$$
\(\sum_{j=0}^\infty |\sa_j(a)(x)|^2 \)^{\f12}=\Bl[\bl(\sum_{2^j\ta
\leq 1}+\sum_{2^j\ta >1}\br) |\sa_j(a)(x)|^2 \Br]^{\f12}\leq
C_{d,p,\phi} \ta^{-d-s} r^{s+d-\f{d-1}p},$$ with  $\ta=d( x,
x_0)$, from which the desired inequality ($\ref{4-2}$) will
follow.

Now the proof of Lemma 4.1 is  reduced to the proof  of the claim
($\ref{4-3}$). We recall that
$$
\sa_j(a)(x)=\int_{\sph} a(y)  G_j(x\cdot y)\, d\sa(y).$$ Hence, by
the definition of regular $(p,\infty,s)$-atom it follows that
\begin{equation}\label{4-4}
|\sa_j(a)(x)|\leq r^{-\f{d-1}p} \int_{B(x_0,r)} \Bl| G_j(x\cdot
y)-\sum_{i=0}^s \f { G_j^{(i)} (x\cdot x_0) } { i!} ( x\cdot
(y-x_0))^i \Br| \, d\sa(y).\end{equation} For $x\in \sph\setminus
B(x_0,4r)$ and $y\in B(x_0,r)$, we write  $\ta =d( x, x_0)$ and
$t=d( x, y)$. Then, evidently, $ \f{3\ta} 4\leq t\leq \f{5\ta}4$
and $ |x\cdot (y-x_0)| =2\bl|\sin\f{\ta-t} 2\sin\f{\ta+t}
2\br|\leq \f 98\ta r.$ Hence, by Lemma 2.2 it follows that for any
$\ell>0$, $x\in \sph\setminus B(x_0,4r)$ and $y\in B(x_0,r)$
$$\Bl|G_j(x\cdot y)-\sum_{i=0}^s \f
{ G_j^{(i)} (x\cdot x_0) } { i!} ( x\cdot (y-x_0))^i\Br| \leq
C_{p,d,\phi} (\ta r)^{s+1} 2^{ j( d+2s+1)} \min\{ 1,
(2^j\ta)^{-\ell} \},$$ where $\ta=d( x, x_0)$. Substituting this
last estimate into ($\ref{4-4}$), we obtain that for $x\in
\sph\setminus B(x_0,4r)$
\begin{align*} |\sa_j(a)(x) |&\le C_{p,d,\phi} \int_{B(x_0, r)} r^{-\f
{d-1}p} (\ta r)^{s+1} 2^{j(d+2s+1)} \min\{
1, (2^j\ta)^{-\ell}\}\, d\sa(y)\\
&\leq  C_{p,d,\phi}2^{j(d+2s+1)} \ta^{s+1} r^{s-\f{d-1}p+d} \min\{
1, (2^j\ta)^{-\ell}\},\end{align*}  proving the claim
($\ref{4-3}$).  This completes the proof.
\end{proof}

{\it Proof of Theorem 1.1 for $0<p\leq 1$.}\    \  Following the
proof in the last section, we need only verify the following three
assertions  in the case $0<p\leq 1$:

(a)\   \  For any sequence $\{a_B\}_{B\in\mathcal{B}}$ of complex
numbers,
\begin{equation*}\Bl\|\(\sum_{B\in\mathcal{B}} |a_B|^2 |\psi_B|^2\)^{\f12}
\Br\|_p \leq C_{p,d}\Bl\|\(\sum_{B\in\mathcal{B}} |a_B|^2|B|^{-1}
\chi_B\)^{\f12} \Br\|_p;\end{equation*}

(b)\    \  for $f\in H^p$,
\begin{equation*}
\Bl\|\(\sum_{B\in\mathcal{B}} \bl|\lb f,\psi_B\rb\br|^2|B|^{-1}
\chi_B\)^{\f12}\Br\|_p\leq C_{p,d} \|f\|_{H^p};\end{equation*}

(c)\    \  for any finite subset $\mathcal{F}\subset \mathcal{B}$
and any sequence $\{a_B\}_{B\in\mathcal{F}} $ of complex numbers,
$$\Bl \|\sum_{B\in\mathcal{F}} a_B \psi_B\Br\|_{H^p}\leq C_{p,d}
\Bl\|(\sum_{B\in\mathcal{F}} |a_B|^2 |\psi_B|^2)^{\f12}\Br\|_p.$$

  Assertion (a) follows directly from the Fefferman-Stein
  inequality and the second inequality in Lemma 3.1 with $0<r<p$.

  For the proof of assertion (b), we define for $\b>0$ and $j\in\mathbb{Z}_+$,
$$ \sa_{j,\b} ^{\ast\ast} (f)(x):=\sup_{ y\in\sph}
\f{|\sa_j(f)(y)|} { \( 1+2^j d( x, y)\)^{\b (d-1)} }.$$ Then by
Lemma 2.3  applied to $N=2^j$,  we have
$$ \sa_{j,\b}^{\ast\ast} (f)(x) \leq C_{d,\b} \( M
(|\sa_j(f)|^{\f1\b})(x) \)^\b,$$ and hence, by Lemma 4.1 and  the
Fefferman-Stein inequality, we obtain that for $\b>\f1p$
\begin{equation}\label{4-5}\Bl\|\(\sum_{j=0}^\infty |\sa_{j,\b}^{\ast\ast} (f)|^2
\)^{\f12} \Br\|_p\leq C_{p,d,\b} \Bl\|\(\sum_{j=0}^\infty |\sa_{j}
(f)|^2 \)^{\f12} \Br\|_p \leq
C_{p,d,\b}'\|f\|_{H^p}.\end{equation} Now assertion (b) follows
from ($\ref{4-5}$) and the following inequality, which can be
easily verified:
$$
\(\sum_{B\in\mathcal{B}}  \bl|\lb f,\psi_B\rb\br|^2|B|^{-1}
\chi_B(x)\)^{\f12}\leq C_{\b,d}  \(\sum_{j=0}^\infty
|\sa_{j,\b}^{\ast\ast} (f) (x) |^2 \)^{\f12},\      \     \forall
\b>0, \  \forall x\in\sph.$$

It remains to prove assertion (c). For simplicity, we define
$a_B=0$ for $B\in\mathcal{B}\setminus \mathcal{F}$,  we denote
 by $\mathcal{D}$ the set of all
spherical caps $B(x_{j,k}, 2^{-j-1})$, ($j\in\mathbb{Z}_+$,
$k\in\Ld_j^d$),  for each $j\in\mathbb{Z}_+$, $k\in\Ld_j^d$, we
write
$$a_{j,k}=a_{B(x_{j,k}, 2^{-j-1})}:= a_{B(x_{j,k}, 2^{-j}\pi)},\
\  \ld_{B(x_{j,k}, 2^{-j-1})}:=\ld_{j,k},\      \ \psi_{B(x_{j,k},
2^{-j-1})}: =\psi_{j,k},$$ and for $c>0$, $B=B(x,r)$, we write
$cB:=B(x,cr)$.

First, we observe that for $0<p\leq 1$
$$\Bl\|\sum_{\sub{2^j\leq 8(d-1) (\f1p-1)\\   k\in\Ld_j^d}}
a_{j,k}\psi_{j,k}\Br\|_{H^p}\leq C_{p,d}\Bl\|\(\sum_{\sub{2^j\leq
8(d-1) (\f1p-1)\\   k\in\Ld_j^d}}
|a_{j,k}|^2|\psi_{j,k}|^2\)^{\f12}\Br\|_{p}$$
 since any two quasi-norms on a finite--dimensional linear space are equivalent.    Thus, without loss of
generality, we may assume
\begin{equation}\label{4-6}
a_{j,k}=0,\     \  \text{for all $2^j\leq 8(d-1)(\f1p-1)$ and
$k\in\Ld_j^d$}.\end{equation}

 Set
$ W(x):= \(\dsum_{B\in \mathcal{D}} |a_B|^2 \ld_B |B|^{-2}
\chi_B(x) \)^{\f12}.$ Then by the first inequality in Lemma 3.1,
$W(x)\leq C_d \(\dsum_{B\in \mathcal{B}} |a_B|^2
|\psi_B(x)|^2\)^{\f12}.$ Thus, it suffices to prove
\begin{equation}\label{4-7}
\|\sum_{B\in\mathcal{B}} a_B \psi_B\|_{H^p} \leq C_{p,d}
\|W\|_p.\end{equation} For the proof of this last inequality, for
an integer  $k\in\mathbb{Z}$ we let
$$ \Og_k=\Bl\{ x\in\sph:\   \   W(x) \ge 2^k\Br\},\      \  \mathcal{A}_k =\Bl\{ B\in\mathcal{D}:\    \   |B\cap
\Og_k|\ge \f12 |B|\Br\}$$ and let $\mathcal{D}_k
=\mathcal{A}_k\setminus \mathcal{A}_{k+1}$. Then  summation by
parts yields \begin{equation}\label{4-8} \(\sum_{k\in\mathbb{Z}}
2^{kp} |\Og_k|\)^{\f1p} \leq C_p \|W\|_p.\end{equation} For each
integer $k\in\mathbb{Z}$,  we choose  a subset
 $\{ B_k^i:\ \ i\in\Ga_k\}$ of $\mathcal{D}_k$ such that
$B_k^i\cap B_k^j=\emptyset$ if $i,j\in\Ga_k$ and $ i\neq j$,  and
$\displaystyle \bigcup_{{B\in\mathcal{D}_k}} B \subset \bigcup_{
i\in\Ga_k} 2B_k^i.$ (The existence of such a subset  $\{ B_k^i:\ \
i\in\Ga_k\}$ is easy to verify.) Now for $k\in\mathbb{Z}$ and
$i\in\Ga_k$, we define
$$ b(k,i) :=|B_k^i|^{\f1p-\f12} \( \sum_{\sub{
B\in\mathcal{D}_k\\
B\subset 2B_k^i}} \ld_B |B|^{-1} |a_B|^2\)^{\f12}$$ and
$$m_{k,i}(x) :=\begin{cases}
b(k,i)^{-1} \dsum_{\sub{ B\in\mathcal{D}_k\\   B\subset 2B_k^i}}
a_B
\psi_B (x),&\      \     \text{if $b(k,i)\neq 0$,}\\
0,&\    \   \text{if $b(k,i)=0$.}\end{cases}$$ Then it's easily
seen that $ \dsum_{B\in\mathcal{D}} a_B\psi_B
=\sum_{k\in\mathbb{Z}}\sum_{i\in\Ga_k} b(k,i) m_{k,i}.$ Therefore,
for the proof of ($\ref{4-7}$), it is sufficient  to show  that
each $m_{k,i}$ is a $p$-molecule up to an absolute constant  and
 \begin{equation}\label{4-9} \(\sum_{k\in\mathbb{Z}}\sum_{i\in\Ga_k}
|b(k,i)|^p\)^{\f1p}\leq C_{p,d} \|W\|_p.\end{equation}

 Since the spaces
$\{\mathcal{H}_k\}_{k=0}^\infty$ (of spherical harmonics) are
mutually orthogonal, it follows  by the assumption ($\ref{4-6}$)
that
\begin{equation*}\label{4-10}
\int_{\sph} m_{k,i} (x) p(x) \, d\sa(x) =0,\    \ \forall  p\in
\Pi_{[(d-1)(\f1p-1)]}.\end{equation*} Thus, to show that $m_{k,i}$
is a $p$-molecule  up to an absolute constant, it suffices to
prove that for any $s>0$,
\begin{equation}\label{4-11} \(
\int_{\sph} |m_{k,i}(x)|^2 \( 1+ \f {d( x, z)} r\)^s\, dx
\)^{\f12} \leq C_{p,d}  r^{(d-1)(\f12-\f1p)},\end{equation} where
$z\equiv z_{k,i}$ denotes  the center of $B_k^i$, $r\equiv
r_{k,i}$ denotes the radius of $B_k^i$. In fact, since
$$ b(k,i)\|m_{k,i}\|_2\leq C_{p,d} \(\sum_{\sub{ B\in\mathcal{D}_k\\
B\subset 2B_k^i}} |a_B|^2 \ld_B |B|^{-1} \)^{\f12}= b(k,i)
|B_k^i|^{\f12-\f1p},$$ it follows that
\begin{equation}\label{4-12}\|m_{k,i}\|_2\leq C_{p,d} r^{(d-1)(\f12-\f1p)}.\end{equation}
 However, on the other hand,  for each $B\in\mathcal{D}_k$ with $B\subset 2 B_k^i$, and
each $x\in\sph\setminus 4B_k^i$,  applying Lemma 2.2, we obtain
$$|\psi_B(x) \leq C_{d,\ell} |B|^{\f\ell{d-1}-1}\sqrt{\ld_B} (d(
x,z))^{-\ell},\    \   \forall\ \ell>\f{s+d-1}2.$$  Hence,
\begin{align*}
\Bl|\sum_{\sub{B\in\mathcal{D}_k\\  B\subset 2B_k^i}} a_B \psi_B
(x)\Br|^2&\leq C_{d,\ell}\( \sum_{\sub{B\in\mathcal{D}_k\\
B\subset
2B_k^i}} \ld_B |B|^{-1} |a_B|^2 \) \(\sum_{\sub{B\in\mathcal{D}_k\\
B\subset 2B_k^i}}
|B|^{-1+\f{2\ell}{d-1}}\)  \( d( x, z)\)^{-2\ell}\\
&\leq C_{d,\ell}  (b(k,i))^2   r^{2\ell-\f{2(d-1)}p}\( d(
x,z)\)^{-2\ell},\    \  \forall\ \ell>\f{s+d-1}2,\end{align*}
which implies
$$ |m_{k,i}(x)|^2\leq C_{d,\ell}     r^{2\ell-\f{2(d-1)}p}\( d(x, z)\)^{-2\ell},\    \  \forall\ \ell>\f{s+d-1}2.$$ We then
deduce by a straightforward calculation that
 $$
\(\int_{d( x, z) \ge 4r} |m_{k,i}(x)|^2 \( 1+\f {d( x,  z)}
r\)^s\, d\sa(x) \)^{\f12}\leq C_{p,d}  r^{(d-1) (\f12-\f1p)},$$
which together with ($\ref{4-12}$) implies ($\ref{4-11}$).

It remains to prove  ($\ref{4-9}$). We observe that
\begin{align*}
\sum_{\sub{B\in\mathcal{D}_k\\
  B\subset 2B_k^i}}\ld_B |B|^{-1}  |a_B|^2 &\leq 2
\sum_{\sub{B\in\mathcal{D}_k\\  B\subset 2B_k^i}} |a_B|^2 |B|^{-2}
\ld_B |B\cap \Og_{k+1}^c| \\
& \leq 2\int_{ (2B_k^i)\cap \Og_{k+1}^c} (W(x))^2 \, d\sa(x) \leq
2^{2k+d+2} |B_k^i|.\end{align*} Thus, by the definition,
$$b(k,i) \leq 2^{(d+2)/2}2^k |B_k^i|^{\f1p},\   \   (k\in\mathbb{Z},\    i\in\Ga_k),$$
which implies
$$\sum_{k,i} |b(k,i)|^p\leq 2^{(d+2)p/2} \sum_{k\in\mathbb{Z}} 2^{kp}
\sum_{ i\in\Ga_k} |B_k^i| \leq 2^{(d+2)p/2} \sum_{k\in\mathbb{Z}}
2^{kp}|\Og_k|$$ since $\{B_k^i\}_{i\in\Ga_k} $ is a sequence of
mutually disjoint subsets of $\Og_k$. This combined with
($\ref{4-8}$) gives ($\ref{4-9}$) and therefore completes the
proof.\hb

\section{Proof of Theorem 1.2}

Recall that $\sa_j$ is defined by ($\ref{2-1}$). We need the
following
\begin{lem}
For $\al>0$ and $ 0<p,\tau\leq \infty$,
$$
|f|_{B_\tau^\al(H^p)}\asymp \( \sum_{j=0}^\infty 2^{j\al\tau}
\|\sa_{j}\circ \sa_j (f)\|^\tau _{H^p} \)^{\f1\tau}\asymp \(
\sum_{j=0}^\infty  2^{j\al\tau} \|\sa_{j} (f)\|^\tau_{H^p}
\)^{\f1\tau},$$ with the usual change when $\tau=\infty$, where
the constants of equivalence depend only on $p, d,\al,\tau$ and
$\phi$.\end{lem}
\begin{proof}
By the definition,  it's easily seen that the series
$\dsum_{j=0}^\infty \sa_{j}\circ\sa_j(f)$ converges to $f$ in the
space $H^p$ and for each $k\in \mathbb{Z}_+$, $\dsum_{j=0}^k
\sa_{j}\circ\sa_j(f)\in \Pi_{2^k-1}$. Thus, for each $k\in
\mathbb{Z}_+$ and $f\in H^p$,
\begin{equation*}
E_{2^k-1} (f)_{H^p} \leq \|\sum_{j=k+1}^\infty \sa_j\circ \sa_j
(f)\|_{H^p}\leq \(\sum_{j=k+1}^\infty \|\sa_j\circ \sa_j
(f)\|_{H^p}^q\)^{\f1q},\end{equation*} where $q=\dmin\{p,1\}.$
Since the operators $\sa_j$, $j\in\mathbb{Z}_+$ are uniformly
bounded on $H^p$,  by the definition and Hardy-type inequality it
follows that
\begin{equation*}
|f|_{B_\tau^\al (H^p)}\leq C_{\tau,\al}  \(\sum_{j=0}^\infty
2^{j\al\tau } \|\sa_j\circ \sa_j(f)\|_{H^p}^\tau \)^{\f1\tau}\leq
C_{p,\tau,\al,\phi} \(\sum_{j=0}^\infty  2^{j\al\tau} \|
\sa_j(f)\|_{H^p}^\tau \)^{\f1\tau},\end{equation*} with the usual
change when $\tau=\infty$.

On the other hand, noting that $\sa_j(g)=0$ for any
$g\in\Pi_{[2^{j-2}]}$ and  $j\ge 1$, we obtain that for $j\ge 1$
$$  \|\sa_j(f)\|_{H^p}
=\inf_{g\in\Pi_{[2^{j-2}]}} \|\sa_j (f-g)\|_{H^p}\leq C_{p,\phi}
E_{[2^{j-2}]}(f)_{H^p}.$$ This together with the uniform
boundedness of the operators $\sa_j$ on $H^p$  implies the desired
inverse  inequalities
$$ \(\sum_{j=0}^\infty  2^{j\al\tau} \|\sa_j\circ \sa_j
(f)\|_{H^p} ^\tau \)^{\f1\tau}\leq C_{p,\phi}
 \(\sum_{j=0}^\infty  2^{j\al\tau} \| \sa_j
(f)\|_{H^p} ^\tau \)^{\f1\tau} \leq C_{p,\phi}'|f|_{B_\tau^\al
(H^p)},$$ with the usual change when $\tau=\infty$. This completes
the proof.
\end{proof}

As indicated in Section 1, the MZ-type inequality ($\ref{1-4}$) in
\mbox{Theorem A} under  the restriction  $0\leq t \leq
\min\{p,1\}$ was proved in [BD, Theorem 3.1].  Our Lemma 5.2 below
asserts that this same inequality, in fact,   holds for
$t>\min\{p,1\}$ as well.

\begin{lem}  Suppose that $\{\xi_k\}_{k\in\Og} $ is a finite
subset of $\sph$ and $\{ a_{N,k}\}_{k\in\Og} $ is a sequence of
nonnegative  numbers smaller than $C_d N^{-(d-1)}$.  Suppose
further that the MZ inequality ($\ref{1-4}$) holds for all
$f\in\Pi_N$ and   $0\leq t \leq \min\{p,1\}$. Then we have, for
any $0<p\leq \infty$, any $f\in \Pi_N$ and all $\al\ge 0$,
$$ \|f\|_p\asymp \begin{cases}
\( \f 1{N^{d-1}} \dsum_{k\in\Og} (N^{d-1} a_{N,k})^{\al p}
|f(\xi_k)|^p\)^{\f1p}, &\      \  \text{if $0<p<\infty$,}\\
\dmax _{k\in\Og} \Bl[ (N^{d-1} a_{N,k})^\al |f(\xi_k)|\Br],&\    \
\text{if $p=\infty$,}
\end{cases}$$
with the constants of equivalence depending only on $d$, $p$ when
$p$ is small, and  $\al$ when $\al$ is  big.\end{lem}

\begin{proof}
Since $0\leq a_{N,k}\leq C_d N^{-(d-1)}$, by ($\ref{1-4}$) with
$t=0$ it follows that
$$ \( \f1{N^{d-1}}\sum_{ k\in\Og} ( N^{d-1} a_{N,k})^{\al p}
|f(\xi_k)|^p \)^{\f1p}\leq C_{p,d,\al} \|f\|_p,\      \  \forall
\al\ge 0,$$ with the usual change when $p=\infty$.  To prove the
inverse inequality, without loss of generality, we may assume
$\al>1$.   Set $\al_0=\min\{p,1\}$. If $0<p<\infty$ then using
($\ref{1-4}$) and H\"older's inequality, we have
\begin{align*}
\|f\|_p^p&\leq C_{p,d}\( \f1{N^{d-1}} \sum_{k\in\Og} |f(\xi_k)|^p
(N^{d-1} a_{N,k})^{\al_0}\) \\
&\leq C_{p,d}\(\f1{N^{d-1}} \sum_{k\in\Og} |f(\xi_k)|^p (N^{d-1}
a_{N,k})^{\al p}\)^{\f{\al_0}{\al p}} \( \f1{N^{d-1}}
\sum_{k\in\Og} |f(\xi_k)|^p\)^{ 1-\f{\al_0}{ \al p}}\\
&\leq C_{p,d}\(\f1{N^{d-1}} \sum_{k\in\Og} |f(\xi_k)|^p (N^{d-1}
a_{N,k})^{\al p}\)^{\f1p\cdot
\f{\al_0}{\al}}\|f\|_p^{p-\f{\al_0}\al};\end{align*} and if
$p=\infty$, using ($\ref{1-4}$) with $t=1$, we have
$$\|f\|_\infty \leq C_d \max_{k\in\Og} \Bl[( N^{d-1}
a_{N,k})|f(\xi_{N,k})|\Br]\leq C_d \(\max_{k\in\Og}( N^{d-1}
a_{N,k})^\al |f(\xi_{N,k})|\)^{\f1\al}\|f\|_\infty ^{1-\f1\al}.$$
Therefore, in either case, we have the desired inverse inequality
$$  \|f\|_p\leq C_{p,d,\al}\( \f1{N^{d-1}}\sum_{ k\in\Og} ( N^{d-1} a_{N,k})^{\al p}
|f(\xi_k)|^p \)^{\f1p},\      \  \forall \al\ge 0,$$ with the
usual change when $p=\infty$.

\end{proof}

{\it Proof of Theorem 1.2.}\    \  We start with the proof of the
equivalence ($\ref{1-10}$). We first note that for each $j\ge 0$,
$\sa_{j}(f) \in \Pi_{2^j}$ and $ \lb f,\psi_{j,k}\rb =
\sqrt{\ld_{j,k}} \sa_j(f)(x_{j,k}).$  Thus, by ($\ref{1-6}$) and
\mbox{Lemma 5.2} with $\al=\f12$ it follows that for $j\ge 0$ and
$0<p\leq \infty$,
\begin{equation}\label{5-1}\|\sa_j(f)\|_p \asymp 2^{-j(d-1)(\f1p-\f12)} \( \dsum_{k\in\Ld_j^d} |\langle
f,\psi_{j,k}\rangle |^p\)^{\f1p},\end{equation} with the usual
change when $p=\infty$.  This together with Lemma 5.1 implies
($\ref{1-10}$) for $1<p\leq \infty$. To show ($\ref{1-10}$) for $
0<p\leq 1$, by Lemma 5.1 it suffices to prove that for $f\in H^p$
and $j\ge 1$
\begin{equation}\label{5-2}
\|\sa_j\circ \sa_j(f)\|_{H^p}\leq
C_{p,d,\phi}\(\sum_{k\in\Ld_j^d}|\langle f,\psi_{j,k}\rangle|^p
2^{j(d-1)(\f p2-1)} \)^{\f1p}\leq C_{p,d,\phi}'
\|\sa_j(f)\|_{H^p}.
\end{equation}
We note that
$$ \sa_j\circ \sa_j (f)(x)=\int_{\sph} \sa_j(f) (y) G_j(x\cdot y)\,
d\sa(y),\    \  x\in\sph.$$ Hence,  by the cubature formula
($\ref{1-5}$) it follows that
$$ \sa_j\circ \sa_j (f)(x)=\sum_{k\in\Ld_j^d} \ld_{j,k} \sa_j(f)(x_{j,k})
G_j(x\cdot x_{j,k})=\sum_{k\in\Ld_j^d}\langle f,\psi_{j,k}\rangle
\psi_{j,k}(x), $$  which together with Theorem 1.1 implies
\begin{equation}
\label{5-3}\|\sa_j\circ \sa_j(f)\|_{H^p} \le
C_{p,d,\phi}'\(\sum_{k\in\Ld_j^d}|\langle f,\psi_{j,k}\rangle|^p
2^{j(d-1)(\f p2-1)} \)^{\f1p},
\end{equation}
since $\dsum_{k\in\Ld_j^d} \chi_{_{B(x_{j,k}, 2^{-j}\pi)}} (x)
\leq C_d$. Thus, combining ($\ref{5-3}$) with ($\ref{5-1}$),
taking into account the fact that $\|\sa_j(f)\|_p\leq
\|\sa_j(f)\|_{H^p}$, we deduce ($\ref{5-2}$).

It remains to prove the inequality ($\ref{1-11}$). Without loss of
generality, we may assume that only a finite number of the
coefficients $a_{j,k}$ are nonzero.  We then deduce from the
definition and the Hardy--type inequality that
\begin{equation}\label{5-4}
|f|_{B_\tau^\al (H^p)} \leq C_{p,\tau}  \( \sum_{j=0}^\infty
2^{j\al\tau} \|\sum_{ k\in \Ld_j^d} a_{j,k}
\psi_{j,k}\|_{H^p}^\tau \)^{\f1\tau},\end{equation} with the usual
change when $\tau=\infty$. However, for $0<p<\infty$,  by Theorem
1.1 we have
\begin{equation}\label{5-5}\|\sum_{k\in\Ld_j^d} a_{j,k}
\psi_{j,k}\|_{H^p}\leq C_{p,d} 2^{j(d-1) (\f12-\f1p)}
\(\sum_{k\in\Ld_j^d} |a_{j,k}|^p\)^{\f1p},
\end{equation} since $\dsum_{k\in\Ld_j^d}\chi_{ _{B(x_{j,k},
2^{-j}\pi)}} \leq  C_d$;  while for  $p=\infty$, we apply Lemma
2.2 to obtain
$$ |\psi_{j,k} (x)|\leq C_{d,\phi}  2^{\f{j(d-1)}2} \min\{ 1, (2^j
d( x, x_{j,k}))^{-d} \},$$ from which it follows that
\begin{align}
|\sum_{k\in\Ld_j^d} a_{j,k} \psi_{j,k}(x)|\leq C_{d,\phi} &
2^{j(d-1)/2} \(\max_{k\in\Ld_j^d} |a_{j,k}| \) \Bl|
\sum_{k\in\Ld_j^d }
 \min\{ 1, (2^j d(x, x_{j,k}))^{-d} \}\Br|\notag  \\
\leq C_{d,\phi} & 2^{j(d-1)/2} \(\max_{k\in\Ld_j^d} |a_{j,k}| \)
\( \sum_{\sub{ x_{j,k}\in B(x, \f \pi {2^j})\\ k\in\Ld_j^d}} 1 +
\sum_{i=1}^{2^j-1} \sum_{\sub{ \f { i\pi} { 2^j} <d( x,
x_{j,k}) \leq \f{ (i+1)\pi}{2^j}\\ k\in \Ld_j^d}} i^{-d}\)\notag  \\
\leq C_{d,\phi} & 2^{j(d-1)/2} \(\max_{k\in\Ld_j^d} |a_{j,k}| \)
\( 1+\sum_{i=1}^\infty i^{-2}\) \leq C_{d,\phi}  2^{j(d-1)/2}
\(\max_{k\in\Ld_j^d} |a_{j,k}| \).\label{5-6}\end{align}  Now
substituting ($\ref{5-5}$) and ($\ref{5-6}$) into ($\ref{5-4}$)
gives the desired inequality ($\ref{1-11}$). This completes the
proof. \hb

\section{Proofs of Theorems 1.3 and 1.4}

{\it Proof of Theorem 1.3.}\    \  For a spherical cap
$B=B(x_{j,k}, 2^{-j}\pi)\in\mathcal{B}$, we set $ w_B(x)=( 1+ 2^j
d( x,  x_{j,k}))^{-d-1}$.
 We define
$$F(x):=\(\sum _{B\in\mathcal{B}} |\langle f,\psi_B\rb |^\tau
|B|^{(-\f12-\f\al{d-1})\tau}w_B(x)\)^{\f1\tau}.$$ Then by Theorem
1.2, $\|F\|_\tau \asymp |f|_{B_\tau^\al(H^\tau)}.$ Thus, it will
suffice to prove
\begin{equation}\label{6-1}
\|f-G_n^p(f)\|_{H^p}\leq C_{p,\al}n^{-\f\al{d-1}}
\|F\|_\tau.\end{equation}

Recall that $$G_n^p(f)=\sum_{B\in\Ga_n} \lb f,\psi_B\rb \psi_B,$$
where $\Ga_n\equiv \Ga_{n, p,f}$ is a set of $n$ spherical caps
$B\in\mathcal{B}$ such that
$$ \min _{B\in\Ga_n} |\lb f,\psi_B\rb| |B|^{\f 1p-\f12}\ge \max _{B\in\mathcal{B}\setminus\Ga_n}
 |\lb f,\psi_B\rb| |B|^{\f 1p-\f12}.$$
 We set
 $$\Sigma(x):=\(\sum_{B\in\mathcal{B}\setminus \Ga_n}|\lb
f,\psi_B\rb|^q|\psi_B(x)|^q\)^{\f1q},$$ where and throughout the
proof, $q=2$ if $p<\infty$, and $=1$ if $p=\infty$. It then
follows by \mbox{Theorem 1.1} that
\begin{equation}\label{6-2}
\|f-G_n^p(f)\|_{H^p} \leq C_{p} \|\Sigma\|_p.
\end{equation} We claim that for $x\in\sph$
\begin{equation}\label{6-3} \Sigma(x) \leq C_{p,\al} n^{-\f\al{d-1}} \bl(F(x)\br)^{\f\tau
p} \|F\|_\tau^{1-\f\tau p},\end{equation} which combined with
($\ref{6-2}$) will  give the desired inequality ($\ref{6-1}$).

To prove ($\ref{6-3}$),  we write, for $t>0$,
$$\Sigma^t_1(x):=\(\sum_{ \sub{B\in\mathcal{B}\setminus \Ga_n\\   |B|\ge t^{d-1}}} |\lb
f,\psi_B\rb|^q|\psi_B(x)|^q\)^{\f1q},$$
$$\Sigma^t_2(x):=\(\sum_{ \sub{B\in\mathcal{B}\setminus \Ga_n\\   |B|< t^{d-1}}} |\lb
f,\psi_B\rb|^q |\psi_B(x)|^q\)^{\f1q}.$$
 We note that for
$B\in \mathcal{B}\setminus \Ga_n$,
\begin{equation}
\label{6-4}|\lb f,\psi_B\rb| |B|^{\f 1p-\f12} \leq  n^{ -\f1\tau}
\( \sum_{ \wt{B}\in \Ga_n}|\lb f,\psi_{\wt{B}}\rb|^\tau
|\wt{B}|^{\f \tau p-\f\tau 2}\)^{\f1\tau} \leq C_{p,\al}n^{
-\f1\tau} \|F\|_{\tau},\end{equation} and by Lemma 2.2,
\begin{equation}\label{6-5}|\psi_B(x)|\leq C_d  |B|^{-\f12}
(w_B(x))^{1+\f1\tau}.\end{equation} Therefore, for
$\Sigma_1^t(x)$, using ($\ref{6-4}$) and ($\ref{6-5}$), we have
\begin{equation}\label{6-6}
\Sigma^t_1(x) \leq C_{p,\al}   \( \sum_{
\sub{B\in\mathcal{B}\setminus\Ga_n \\ |B|\ge t^{d-1}}} |B|^{-\f
qp} (w_B(x))^{(1+\f1\tau)q}\)^{\f1q} n^{-\f 1\tau}\|F\|_\tau \leq
C_{p,\al} n^{-\f 1\tau} t^{-\f{d-1}p} \|F\|_\tau,\end{equation}
while for $\Sigma^t_2(x)$, using ($\ref{6-5}$),  we have
\begin{align}
\Sigma^t _2(x)&\leq C_{p,\al}\(\sum_{ \sub{B\in\mathcal{B}\setminus\Ga_n\notag\\
|B|\le t^{d-1}}}|\lb f,\psi_B\rb|^q|B|^{-\f q2-\f{q\al}{d-1}}
|B|^{\f{q\al}{d-1}}( w_B(x))^{(1+\f1\tau)q} \)^{\f1q} \notag \\
&\leq C_{p,\al} \max_{B\in\mathcal{B}} \(|\lb
f,\psi_B\rb||B|^{-\f12-\f{\al}{d-1}}(w_B(x))^{\f1\tau}\)
\(\sum_{\sub{
B\in\mathcal{B}\notag\\
|B|\le t^{d-1}}} |B|^{\f{q\al}{d-1}} (w_B(x))^q\)^{\f1q}\\
&\le C_{p,\al}t^\al |F(x)|.\label{6-7}\end{align} Combining
($\ref{6-6}$) with ($\ref{6-7}$), we obtain $$ \Sigma(x) \leq
C_{p,\al}\Bl[ n^{-\f1\tau} t^{-\f{d-1}p} \|F\|_\tau +t^\al
|F(x)|\Br]$$ and the claim ($\ref{6-3}$) then  follows by taking
$t= n^{-\f1{d-1}} \(\f{\|F\|_\tau}{F(x)}\)^{\f\tau {d-1}}$. This
completes the proof. \hb

{\it Proof of Theorem 1.4.}\    \   \    First, we prove that for
$f=\dsum_{B\in\mathcal{B}} a_B\psi_B$ with only  $n$ nonzero
coefficients $a_B$,
\begin{equation}
\label{6-8}|f|_{B_\tau^\al (H^\tau)} \le C_{p,\al} n^{\f\al{d-1}}
\|g\|_p,\end{equation} where $g(x): =\dmax_{B\in\Ld}
\(|a_B||B|^{-1/2} \    \chi_B(x)\)$  and $\Ld=\{ B\in\mathcal{B}:\
a_B\neq 0\}$.  Suppose $\Ld=\{B_1,\ldots, B_n\}$ with  $ |B_1|\leq
\ldots \leq |B_n|.$ Set $E_1=B_1$ and $\displaystyle E_j =B_j
\setminus \bigcup_{i=1}^{j-1} B_i$, $j\ge 2$.  Then by Theorem
1.2, we have
\begin{align*}
|f|_{B_\tau^\al(H^\tau)} ^\tau &\le C_{p,\al} \int_{\sph}
\(\sum_{j=1}^n |a_{_{B_j}} |^\tau
|B_j|^{-\f\tau2-\f{\al\tau}{d-1}} \chi_{B_j}(x)
\)\,d\sa(x)\\
& =C_{p,\al}\sum_{i=1}^n \int_{E_i} \sum_{j=i}^n |a_{B_j}|^\tau
|B_j|^{-\f\tau2-\f{\al\tau}{d-1}} \chi_{B_j} (x)  \,
d\sa(x)\\
&\leq C_{p,\al} \sum_{i=1}^n \int_{E_i} \(\sum_{\sub{ B\in\mathcal{B}\\
|B|\ge |B_i|}}  |B|^{-\f{\al\tau}{d-1}} \chi_{B} (x) \)
|g(x)|^\tau \,
d\sa(x)\\
&\le C_{p,\al} \sum_{i=1}^n |B_i|^{-\f{\al\tau} {d-1}} \int_{E_i}
|g(x)|^\tau\, d\sa(x)\\
& \le C_{p,\al}  n^{1-\f\tau p} \|g\|_p^\tau= C_{p,\al}
n^{\f{\al\tau}{d-1}} \|g\|_p^\tau,
\end{align*}
which gives  ($\ref{6-8}$).

Next, we show
$$ \|G_n^p(f)\|_{_{B_\tau^\al (H^\tau)}} \le C_{p,\al} n^{\f\al{d-1}}
\|f\|_{H^p}.$$  By ($\ref{6-8}$), it suffices to prove that
\begin{equation}
\label{6-9}\Bl\|\dmax_{B\in\Ga_n} \(|\lb f,\psi_B\rb||B|^{-1/2}
\chi_B(x)\)\Br\|_p\leq C_{p,\al} \|f\|_{H^p}.\end{equation}
($\ref{6-9}$) for $0<p<\infty$  follows directly from Theorem 1.1,
while for $p=\infty$ follows by  the following inequality, which
can be easily deduced by  Lemma 2.2:
$$\max_{B\in\mathcal{B}} |\lb f,\psi_B\rb | |B|^{-\f12} \chi_B(x)
\leq C_{d,\phi} M(f) (x),\      \    x\in\sph.$$ This completes
the proof.\hb
\section*{Acknowledgements}
The author would like to express his sincere gratitude to an
anonymous referee for many helpful comments on this paper.




\end{document}